\documentclass[10pt]{article}
\usepackage{amsfonts}
\usepackage{hyperref}
\usepackage{mathrsfs,enumerate,amssymb,amsmath,color,lineno}
\usepackage{indentfirst}
\usepackage{amsmath}
\usepackage{amssymb}
\usepackage[amsmath,thmmarks]{ntheorem}
\usepackage{graphicx}
\usepackage{tikz}
\usepackage{varioref}

\newtheorem{definition}{\bf Definition}[section]
\newtheorem{lemma}{\bf Lemma}[section]
\newtheorem{theorem}{\bf Theorem}[section]
\newtheorem{remark}{\bf Remark}[section]
\newtheorem{corollary}{\bf Corollary}[section]
\newtheorem{example}{\bf Example}[section]

\newtheorem{proposition}{\bf Proposition}[section]

\begin{document}
\setcounter{page}{1}

\title{{\textbf{New characterizations of migrative 2-uninorms}}
}
\author{Shudi Liang\footnote{\emph{E-mail address}: shudi$\_$liang@163.com}, Xue-ping Wang\footnote{Corresponding author. xpwang1@hotmail.com; fax: +86-28-84761502}, \\
\emph{School of Mathematical Sciences, Sichuan Normal University,}\\
\emph{Chengdu 610066, Sichuan, People's Republic of China}}

\newcommand{\pp}[2]{\frac{\partial #1}{\partial #2}}
\date{}
\maketitle
\begin{quote}
{\bf Abstract.} This article pays attention to the $\alpha$-migrativity of 2-uninorms with $\alpha\in [0,1]$ deeply. It describes the $\alpha$-migrativity of 2-uninorms completely, which generalizes and unifies some current existing characterizations for $\alpha$-migrativity of triangular norm, triangular conorm, uninorms, nullnorms, uni-nullnorms and null-uninorms, respectively.


\end{quote}

\section{Introduction}
There are many special classes of aggregation operators, such as triangular norms (briefly t-norms), t-conorms, t-operators, nullnorms, uninorms, uni-nullnorms, 2-uninorms, copulas and so on. When using an aggregation technology to solve practical problems, the first problem we face is how to select the most appropriate aggregation function. The selection of aggregation function is determined by its application circumstances, and there is no unified criterions. For example, nullnorms can be applied in the areas of neural networks and decision making \cite{HGC2008,LJ2007}, and uninorms can be used in various fields like neural networks, decision making, fuzzy logics and fuzzy system modeling \cite{TM1998,Y2001,Y2003}. On the other hand, 2-uninorms specifically unify uninorms and nullnorms. They were first investigated by Akella \cite{A2007}, since then, scholars on 2-uninorms have mainly concentrated on the structure of 2-uninorms \cite{ZSLB2018}, distributivity of 2-uninorms \cite{DR2016,SL2020,WQ2017,ZQ2020}, migrativity of 2-uninorms \cite{HQ2021,LW2023,WZZL2019,YQ2022} and $(U^2,N)$-implications of 2-uninorms \cite{ZL2020}, etc.. 2-uninorms are crucial for its practical applications in fuzzy quantifiers or networks as well as theoretical interest \cite{HQ2021,WQ2017,ZL2020} because it is closely related to other aggregation operators.

The study of $\alpha$-migrativity of aggregation function has great significance and value in theory and applications. In applications, the migrativity means that when the input is scaled down during the aggregation process, the inherent properties of the aggregation function itself will not be changed \cite{MBF2010}. So the migrativity of aggregation function is of great importance in image processing and decision-making \cite{BQ2013, BMM2009, BBFMM2012}.
The concept of migrativity comes from an open question, which was proposed by Endre Pap \cite{MN1995,MN1996}: does there exist a strictly increasing t-norm, that is different from the t-norm given by
$$
T(x,y)=
\begin{cases}
\alpha xy, & (x,y)\in (0,1),\\
\min(x,y), & \mbox{otherwise},
\end{cases}
$$
 for the following equation
\begin{equation}\label{(1)}
T(\alpha x,y)=T(x,\alpha y) \mbox{ for any $x,y\in [0,1]$}
\end{equation} for an $\alpha\in (0,1)$? This problem was positively solved by Budin\v{c}evi\'{c} and Kurili\'{c} \cite{BK1998}. In 2008, the property of $T$ satisfying the equation \eqref{(1)} was called an $\alpha$-migrativity of $T$ by Durante and Sarkoci \cite{DS2008}. By replacing the product $\alpha x$ in the equation \eqref{(1)} by $T_1(\alpha, x)$ with $T_1$ a t-norm, Fodor and Rudas \cite{FR2011} generalized the $\alpha$-migrativity to the $(\alpha,T_1)$-migrativity of t-norm $T$, which is expressed as follows
\begin{equation}\label{(2)}
T(T_1(\alpha, x),y)=T(x,T_1(\alpha, y)) \mbox{ for each $x,y\in [0,1]$.}
\end{equation}
Indeed, many scholars all over the world explore the migrative property of several binary aggregation operations and their works can be mainly divided into two categories: one is to change the form of Formula \eqref{(1)} or \eqref{(2)}; the other is to find more general aggregation functions to replace the t-norms $T$ or $T_1$ in Formula \eqref{(2)}. For example, the migrativity of semicopulas, quasicoupulas, copulas \cite{DFQ2012,MBF2010}, fuzzy implications \cite{BJM2020},  semi-t-operators \cite{ZQ2021}, t-subnorms, t-norms, t-conorms \cite{FR2011,O2013,SZQZ2019,WO2013}, uninorms, nullnorms \cite{LQ2021a,LQ2021b, MMRT2013,MMRT2015,SZLZ2015,SZLX2015,SZL2016,SLRRT2017,ZL2017,ZSL2014,ZSL2016}, overlap functions \cite{BDBB2013,QH2018,QH2019,ZH2020,ZY2021,ZWY2021a,ZWY2022,ZZQ2022},
uni-nullnorms \cite{SQZ2019,SQ2021} and 2-uninorms \cite{HQ2021,LW2023,WZZL2019} have been widely studied, respectively. Because t-norm, t-conorm, uninorm, nullnorm, uni-nullnorm and null-uninorm are just subclasses of 2-uninorm and the $\alpha$-migrativity of 2-uninorms is very valuable at least in image processing and decision making \cite{BQ2013, BMM2009, BBFMM2012}, Liang and Wang \cite{LW2023} recently characterized $(\alpha,\mathbb{U}_2)$-migrative 2-uninorms $\mathbb{U}_1$ when $e_1=e_2$ or $a_1=a_2$ or $f_1=f_2$ where $\{e_1,f_1\}_{a_1}$ and $\{e_2,f_2\}_{a_2}$ are 2-neutral elements of $\mathbb{U}_1$ and $\mathbb{U}_2$, respectively. However, the characterization for $(\alpha,\mathbb{U}_2)$-migrative 2-uninorms $\mathbb{U}_1$ when $e_1\neq e_2$ and $a_1\neq a_2$ and $f_1\neq f_2$ has not been resolved yet. Therefore, a complete characterization for $\alpha$-migrativity of 2-uninorms is necessary. This article will focus on this problem.

The rest of this article is organized as follows. In Section 2, we recall some necessary notions and one known result. In Section 3, we describes the migrative 2-uninorms completely. Conclusions
are drawn in Section 4.

\section{Preliminaries}\label{se2}
This section just recall some necessary definitions and results of 2-uninorms. Readers can refer to \cite{KMP2000},\cite{YR1996}, \cite{CB2001} and \cite{SWQ2017} for more information on t-norms, t-conorms, uninorms, nullnorms, uni-nullnorms and null-uninorms, respectively.

%
%
%
%

\begin{definition}[\cite{A2007}]\label{def2.4}
\emph{Let $\mathbb{U}:[0,1]^{2}\to [0,1]$ be a commutative binary operator. Then $\{e,f\}_a$ with $0\leq e\leq a\leq f\leq 1$  is called
a {\it 2-neutral element} of $\mathbb{U}$ if $\mathbb{U}(e,x)=x$ for all $x\in [0, a]$ and $\mathbb{U}(f,x)=x$ for all $x\in [a, 1]$.}
\end{definition}

\begin{definition}[\cite{A2007}]\label{def2.5}
\emph{A binary operation $\mathbb{U}: [0,1]^2 \to [0,1]$ is called a {\it 2-uninorm} if it is commutative, associative, increasing in both variables and has a 2-neutral element $\{e,f\}_a$.}
\end{definition}

%

\begin{definition}[\cite{LW2023}]
\emph{
Let $\mathbb{U}_1$ and $\mathbb{U}_2$ be two 2-uninorms on the unit interval $[0,1]$ and $\alpha \in [0,1]$. Then $\mathbb{U}_1$ is called $(\alpha, \mathbb{U}_2)$-migrative or $\alpha$-migrative over
$\mathbb{U}_2$ if
\begin{align}\label{eq-def}
\mathbb{U}_1(\mathbb{U}_2(\alpha,x),y)=\mathbb{U}_1(x,\mathbb{U}_2(\alpha,y))
\end{align}for any $x,y\in [0,1]$.}
\end{definition}

Let $\lambda=\mathbb{U}_2(\alpha,e_1)$ and $\mu=\mathbb{U}_2(\alpha,f_1)$.
If $\mathbb{U}_1$ is $(\alpha,\mathbb{U}_2)$-migrative, then $\mathbb{U}_2(\alpha,x)=\mathbb{U}_1(x,\lambda)$ for any $\mathbb{U}_2(\alpha,x)\leq a_1$ and  $\mathbb{U}_2(\alpha,x)=
\mathbb{U}_1(x,\mu)$ for any $\mathbb{U}_2(\alpha,x)\geq a_1$, which will be used repeatedly.


\begin{proposition}[\cite{LW2023}]\label{pro-2.3}
Let $\mathbb{U}$ be a 2-uninorm on the unit interval $[0,1]$ with a 2-neutral element $\{e,f\}_a$. Then
\begin{enumerate}[{\rm (i)}]
  \item
$\mathbb{U}(x,y)=\mathbb{U}(x,a)$ for all $x\in [0,a]$ and $y\in [a,f]$.
  \item
$\mathbb{U}(x,y)=\mathbb{U}(a,y)$ for all $x\in [a,f]$ and $y\in [0,a]$.
  \item
$\mathbb{U}(x,y)=\mathbb{U}(x,a)$ for all $x\in [a,1]$ and $y\in [e,a]$.
  \item
$\mathbb{U}(x,y)=\mathbb{U}(a,y)$ for all $x\in [e,a]$ and $y\in [a,1]$.
  \item
$\min(x,y)\leq \mathbb{U}(x,y)\leq a$ for all $(x,y)\in [0,e]\times [e,f]\cup [e,f]\times [0,e]$.
  \item
$a\leq \mathbb{U}(x,y)\leq \max(x,y)$ for all $(x,y)\in [e,f]\times [f,1]\cup [f,1]\times [e,f]$.
 \item
$\mathbb{U}(x,y)\in [0,e)\cup \{a\}\cup (f,1]$ for all $(x,y)\in [0,a]\times [a,1]\cup [a,1]\times [0,a]$.
\end{enumerate}
\end{proposition}

\section{$\alpha$-migrativity of 2-uninorms}\label{sec4}

In this section, we completely characterize the $\alpha$-migrativity of 2-uninorms with 2-neutral elements $\{e_1,f_1\}_{a_1}$ and $\{e_2,f_2\}_{a_2}$, respectively.

\begin{lemma}\label{lemma3.1}
Let $\mathbb{U}_1$ and $\mathbb{U}_2$ be two 2-uninorms on the unit interval $[0,1]$ with 2-neutral elements $\{e_1,f_1\}_{a_1}$ and $\{e_2,f_2\}_{a_2}$, respectively.
\begin{enumerate}[{\rm (i)}]
  \item
If $\mathbb{U}_2(\alpha,e_1)\leq e_1$, then $\mathbb{U}_2(\alpha,f_1)\leq f_1$.
  \item
If $\mathbb{U}_2(\alpha,f_1)\geq f_1$, then $\mathbb{U}_2(\alpha,e_1)\geq e_1$.
\end{enumerate}
\end{lemma}
\begin{proof}
(i) We prove that $\mathbb{U}_2(\alpha,f_1)\leq f_1$ by distinguishing two cases. If $e_1\leq a_2$, then from $\mathbb{U}_2(\alpha,e_1)\leq e_1$ we know $\alpha\leq e_2$, and subsequently, $\mathbb{U}_2(\alpha,f_1)\leq \mathbb{U}_2(e_2,f_1)=f_1$ for each $f_1\leq a_2$ and $\mathbb{U}_2(\alpha,f_1)\leq \mathbb{U}_2(e_2,f_1)\leq \mathbb{U}_2(f_2,f_1)=f_1$ for each $f_1>a_2$, i.e., $\mathbb{U}_2(\alpha,f_1)\leq f_1$.
If $e_1>a_2$, then by $\mathbb{U}_2(\alpha,e_1)\leq e_1$ we get $\alpha\leq f_2$, showing that $\mathbb{U}_2(\alpha,f_1)\leq \mathbb{U}_2(f_2,f_1)=f_1$.

(ii)
We prove that $\mathbb{U}_2(\alpha,e_1)\geq e_1$ by the following two cases. If $f_1\geq a_2$, then from $\mathbb{U}_2(\alpha,f_1)\geq f_1$ we obtain $\alpha\geq f_2$, and subsequently, $\mathbb{U}_2(\alpha,e_1)\leq \mathbb{U}_2(f_2,e_1)=e_1$ for each $e_1\geq a_2$ and $\mathbb{U}_2(\alpha,e_1)\geq \mathbb{U}_2(f_2,e_1)\geq \mathbb{U}_2(e_2,e_1)=e_1$ for each $e_1<a_2$.
If $f_1<a_2$, then by $\mathbb{U}_2(\alpha,f_1)\geq f_1$ we get $\alpha\geq e_2$, showing that $\mathbb{U}_2(\alpha,e_1)\geq \mathbb{U}_2(e_2,e_1)=e_1$.
\end{proof}

\begin{lemma}\label{lemma3.2}
Let $\mathbb{U}_1$ and $\mathbb{U}_2$ be two 2-uninorms on the unit interval $[0,1]$ with 2-neutral elements $\{e_1,f_1\}_{a_1}$ and $\{e_2,f_2\}_{a_2}$, respectively. Suppose that $\mathbb{U}_1$ is $(\alpha,\mathbb{U}_2)$-migrative. The following two statements hold.
\begin{enumerate}[{\rm (i)}]
  \item
 If $\mathbb{U}_2(\alpha,e_1)\in [e_1,a_1]$, then $\mathbb{U}_2(\alpha,a_1)\geq a_1$.
   \item
If $\mathbb{U}_2(\alpha,f_1)\in [a_1, f_1]$, then $\mathbb{U}_2(\alpha,a_1)\leq a_1$.
\end{enumerate}
\end{lemma}
\begin{proof}
Concerning (i), assuming that $\mathbb{U}_1$ is $(\alpha,\mathbb{U}_2)$-migrative. If $\mathbb{U}_2(\alpha,a_1)<a_1$, then
$\mathbb{U}_2(\alpha,a_1)=\mathbb{U}_1(\mathbb{U}_2(\alpha,a_1),e_1)=\mathbb{U}_1(a_1,\mathbb{U}_2(\alpha,e_1))=a_1$, a contradiction. Thus $\mathbb{U}_2(\alpha,a_1)\geq a_1$.

To show (ii), assuming that $\mathbb{U}_1$ is $(\alpha,\mathbb{U}_2)$-migrative. If $\mathbb{U}_2(\alpha,a_1)>a_1$, then $$\mathbb{U}_2(\alpha,a_1)=\mathbb{U}_1(\mathbb{U}_2(\alpha,a_1),f_1)=\mathbb{U}_1(a_1,\mathbb{U}_2(\alpha,f_1))=a_1,$$ a contradiction. Thus $\mathbb{U}_2(\alpha,a_1)\leq a_1$.
\end{proof}

\begin{proposition}\label{pro-3.1}
Let $\mathbb{U}_1$ and $\mathbb{U}_2$ be two 2-uninorms on the unit interval $[0,1]$ with 2-neutral elements $\{e_1,f_1\}_{a_1}$ and $\{e_2,f_2\}_{a_2}$, respectively.
\begin{enumerate}[{\rm (i)}]
  \item
If $\mathbb{U}_2(\alpha,e_1)\geq a_1$, then $\mathbb{U}_1$ is $(\alpha,\mathbb{U}_2)$-migrative iff $\mathbb{U}_1(\mu,x)=\mathbb{U}_2(\alpha,x)$ for each $x\in [0,1]$.
  \item
If $\mathbb{U}_2(\alpha,f_1)\leq a_1$, then $\mathbb{U}_1$ is $(\alpha,\mathbb{U}_2)$-migrative iff $\mathbb{U}_1(\lambda,x)=\mathbb{U}_2(\alpha,x)$ for each $x\in [0,1]$.
\end{enumerate}
\end{proposition}
\begin{proof}

\begin{enumerate}[{\rm (i)}]
\item We first prove $\mathbb{U}_2(\alpha,x)=\mathbb{U}_1(\mu,x)$ for each $x\in [0,1]$ done by checking the subsequent two cases.

Case 1. There is at least one $k\in [0,e_1)$ such that $\mathbb{U}_2(\alpha,k)\leq a_1$.

Let $\mathcal{B}=\{k\in [0,e_1)|\mathbb{U}_2(\alpha,k)\leq a_1\}$ and $\eta=\sup \mathcal{B}$. In this case, the proof is split into two subcases as follows.

Subcase 1. If $\eta\in \mathcal{B}$, then $\mathbb{U}_2(\alpha,x)\leq \mathbb{U}_2(\alpha,\eta)\leq a_1$ for each $x\in [0,\eta]$ and $\mathbb{U}_2(\alpha,x)> a_1$ for each $x\in (\eta,1]$, showing that
$$
\mathbb{U}_2(\alpha,x)=
\begin{cases}
\mathbb{U}_1(\lambda,x), & x\in [0,\eta],\\
\mathbb{U}_1(\mu,x), & x\in (\eta,1].
\end{cases}
$$
Moreover, since $\lambda=\mathbb{U}_2(\alpha,e_1)\geq a_1$, $\mu=\mathbb{U}_2(\alpha,f_1)\geq \mathbb{U}_2(\alpha,e_1)\geq a_1$, by \eqref{eq-def} we have $\mathbb{U}_1(\mu,x)=\mathbb{U}_1(\mathbb{U}_1(\mu,f_1),x)=\mathbb{U}_1(\mathbb{U}_2(\alpha,f_1),x)
=\mathbb{U}_1(f_1,\mathbb{U}_2(\alpha,x))=\mathbb{U}_1(f_1,\mathbb{U}_1(\lambda,x))=
\mathbb{U}_1(\mathbb{U}_1(f_1,\lambda),x)=\mathbb{U}_1(\lambda,x)$ for each $x\in [0,\eta]$ and, consequently, $\mathbb{U}_2(\alpha,x)=\mathbb{U}_1(\mu,x)$ for each $x\in [0,1]$.

Subcase 2. If $\eta\notin \mathcal{B}$, then $\mathbb{U}_2(\alpha,\eta)>a_1$. This follows from the monotonicity of $\mathbb{U}_2$ that $\mathbb{U}_2(\alpha,x)\geq \mathbb{U}_2(\alpha,\eta)>a_1$ for each $x\in [\eta,1]$. Now, let $x\in [0,\eta)$. Then from the definition of $\eta$ there is a $k_0\in \mathcal{B}$ such that $k_0>x$, showing that $\mathbb{U}_2(\alpha,x)\leq \mathbb{U}_2(\alpha,k_0)\leq a_1$ from the monotonicity of $\mathbb{U}_2$. Therefore,
$$
\mathbb{U}_2(\alpha,x)=
\begin{cases}
\mathbb{U}_1(\lambda,x), & x\in [0,\eta),\\
\mathbb{U}_1(\mu,x), & x\in [\eta,1].
\end{cases}
$$
Moreover, similar to Subcase 1 we have $\mathbb{U}_1(\mu,x)=\mathbb{U}_1(\lambda,x)$ for each $x\in [0,\eta)$ and consequently, $\mathbb{U}_2(\alpha,x)=\mathbb{U}_1(\mu,x)$ for all $x\in [0,1]$.

Case 2. There is no $k\in [0,e_1)$ such that $\mathbb{U}_2(\alpha,k)\leq a_1$.

Obviously, $\mathbb{U}_2(\alpha,x)\geq \mathbb{U}_2(\alpha,0)>a_1$ for each $x\in [0,1]$. Thus $\mathbb{U}_2(\alpha,x)=\mathbb{U}_1(\mu,x)$ for each $x\in [0,1]$.

Now, we prove that $\mathbb{U}_1$ is $(\alpha, \mathbb{U}_2)$-migrative. In fact, if $\mathbb{U}_2(\alpha,x)=\mathbb{U}_1(\mu,x)$ for each $x\in [0,1]$, then since $\mathbb{U}_1$ is commutative and associative, we have
$$\mathbb{U}_1(\mathbb{U}_2(\alpha,x),y)=\mathbb{U}_1(\mathbb{U}_1(\mu,x),y)=
\mathbb{U}_1(x,\mathbb{U}_1(\mu,y))=\mathbb{U}_1(x,\mathbb{U}_2(\alpha,y)).$$
  \item
Necessity. This will be done by checking the subsequent two cases.

Case 1. There is at least one $k\in (f_1,1]$ such that $\mathbb{U}_2(\alpha,k)\geq a_1$.

Let $\mathcal{A}=\{k\in (f_1,1]|\mathbb{U}_2(\alpha,k)\geq a_1\}$ and $\xi=\inf \mathcal{A}$. In this case, the proof is split into two subcases as follows.

Subcase 1. If $\xi\in \mathcal{A}$,
then $\mathbb{U}_2(\alpha,x)\geq \mathbb{U}_2(\alpha,\xi)\geq a_1$ for each $x\in [\xi,1]$ and $\mathbb{U}_2(\alpha,x)< a_1$ for $x\in [0,\xi)$, showing that
$$
\mathbb{U}_2(\alpha,x)=
\begin{cases}
\mathbb{U}_1(\lambda,x), & x\in [0,\xi),\\
\mathbb{U}_1(\mu,x), & x\in [\xi,1].
\end{cases}
$$
Moreover, since $\lambda=\mathbb{U}_2(\alpha,e_1)\leq \mathbb{U}_2(\alpha,f_1)\leq a_1$ and $\mu=\mathbb{U}_2(\alpha,f_1)\leq a_1$ we have
\begin{eqnarray*}
\mathbb{U}_1(\lambda,x)&=&\mathbb{U}_1(\mathbb{U}_1(\lambda,e_1),x)\\
&=&\mathbb{U}_1(\mathbb{U}_2(\alpha,e_1),x)\\
&=&\mathbb{U}_1(e_1,\mathbb{U}_2(\alpha,x))\\
&=&\mathbb{U}_1(e_1,\mathbb{U}_1(\mu,x))\\
&=&\mathbb{U}_1(\mathbb{U}_1(e_1,\mu),x)\\
&=&\mathbb{U}_1(\mu,x)
\end{eqnarray*}
for each $x\in [\xi,1]$, and consequently, $\mathbb{U}_2(\alpha,x)=\mathbb{U}_1(\lambda,x)$ for each $x\in [0,1]$.

Subcase 2. If $\xi\notin \mathcal{A}$, then $\mathbb{U}_2(\alpha,\xi)<a_1$. This follows from the monotonicity of $\mathbb{U}_2$ that $\mathbb{U}_2(\alpha,x)\leq \mathbb{U}_2(\alpha,\xi)< a_1$ for each $x\in [0,\xi]$. Now, let $x\in (\xi,1]$. Then from the definition of $\xi$ there exists a $k_0\in \mathcal{A}$ such that $k_0<x$. Thus from the monotonicity of $\mathbb{U}_2$ we have $\mathbb{U}_2(\alpha,x)\geq \mathbb{U}_2(\alpha,k_0)\geq a_1$. Therefore,
$$
\mathbb{U}_2(\alpha,x)=
\begin{cases}
\mathbb{U}_1(\lambda,x), & x\in [0,\xi],\\
\mathbb{U}_1(\mu,x), & x\in (\xi,1].
\end{cases}
$$
Moreover, similar to Subcase 1 we have $\mathbb{U}_1(\lambda,x)=\mathbb{U}_1(\mu,x)$ for each $x\in (\xi,1]$ and consequently, $\mathbb{U}_2(\alpha,x)=\mathbb{U}_1(\lambda,x)$ for each $x\in [0,1]$.

Case 2. There is no $k\in (f_1,1]$ such that $\mathbb{U}_2(\alpha,k)\geq a_1$.

Obviously, $\mathbb{U}_2(\alpha,x)\leq \mathbb{U}_2(\alpha,1)<a_1$ for each $x\in [0,1]$. Therefore, $\mathbb{U}_2(\alpha,x)=\mathbb{U}_1(\lambda,x)$ for each $x\in [0,1]$.

Sufficiency. If $\mathbb{U}_2(\alpha,x)=\mathbb{U}_1(\lambda,x)$ for each $x\in [0,1]$, then since $\mathbb{U}_1$ is commutative and associative, we have
$$\mathbb{U}_1(\mathbb{U}_2(\alpha,x),y)=\mathbb{U}_1(\mathbb{U}_1(\lambda,x),y)=
\mathbb{U}_1(x,\mathbb{U}_1(\lambda,y))=\mathbb{U}_1(x,\mathbb{U}_2(\alpha,y)),$$
i.e., $\mathbb{U}_1$ is $(\alpha, \mathbb{U}_2)$-migrative.
\end{enumerate}
\end{proof}

\begin{proposition}\label{pro-3.3}
Let $\mathbb{U}_1$ and $\mathbb{U}_2$ be two 2-uninorms on the unit interval $[0,1]$ with 2-neutral elements $\{e_1,f_1\}_{a_1}$ and $\{e_2,f_2\}_{a_2}$, respectively.
\begin{enumerate}[{\rm (i)}]
  \item
If $\mathbb{U}_2(\alpha,e_1)\in [e_1,a_1]$, then $\mathbb{U}_1$ is $(\alpha, \mathbb{U}_2)$-migrative iff
\begin{equation*}
\mathbb{U}_2(\alpha,x)=
\begin{cases}
\mathbb{U}_1(\lambda,x), & x\in [0,{a_1}],\\
\mathbb{U}_1(\mu,x), & x\in [{a_1},1]\\
\end{cases}
\end{equation*}
and $\mathbb{U}_1(\mu,x)=\mathbb{U}_1({a_1},x)$ for each $x\in [0,{a_1}]$.
  \item
If $\mathbb{U}_2(\alpha,f_1)\in [a_1,f_1]$, then $\mathbb{U}_1$ is $(\alpha, \mathbb{U}_2)$-migrative iff
\begin{equation*}
\mathbb{U}_2(\alpha,x)=
\begin{cases}
\mathbb{U}_1(\lambda,x), & x\in [0,{a_1}],\\
\mathbb{U}_1(\mu,x), & x\in [{a_1},1]\\
\end{cases}
\end{equation*}
and $\mathbb{U}_1(\lambda,x)=\mathbb{U}_1({a_1},x)$ for each $x\in [a_1,1]$.
\end{enumerate}
\end{proposition}
\begin{proof}
\begin{enumerate}[{\rm (i)}]
  \item
Assume that $\mathbb{U}_1$ is $(\alpha, \mathbb{U}_2)$-migrative. Then by Lemma \ref{lemma3.2} (i) we get $\mathbb{U}_2(\alpha,a_1)\geq a_1$.
Since $\mathbb{U}_2(\alpha,x)\geq \mathbb{U}_2(\alpha,a_1)\geq a_1$ for any $x\geq a_1$, we have $\mathbb{U}_2(\alpha,x)=\mathbb{U}_1(x,\mu)$ and $\mu=\mathbb{U}_2(\alpha,f_1)\geq \mathbb{U}_2(\alpha,a_1)\geq a_1$. From Proposition \ref{pro-2.3} (iii) we get $\mathbb{U}_2(\alpha,x)\leq \mathbb{U}_2(\alpha,a_1)=\mathbb{U}_1(\mu,a_1)=\mathbb{U}_1(\mu,e_1)=
\mathbb{U}_1(\mathbb{U}_2(\alpha,f_1),e_1)
=\mathbb{U}_1(f_1,\mathbb{U}_2(\alpha,e_1))=a_1$ for any $x\leq a_1$. Consequently,
$$
\mathbb{U}_2(\alpha,x)=
\begin{cases}
\mathbb{U}_1(\lambda,x), & x\in [0,a_1],\\
\mathbb{U}_1(\mu,x), & x\in [a_1,1].\\
\end{cases}
$$

Next, from \eqref{eq-def} we have
\begin{eqnarray*}
\mathbb{U}_1(\mu,x)&=&\mathbb{U}_1(\mathbb{U}_1(\mu,f_1),x)\\&=&\mathbb{U}_1(\mathbb{U}_2(\alpha,f_1)
,x)\\&=&\mathbb{U}_1(f_1,\mathbb{U}_2(\alpha,x))\\&=&\mathbb{U}_1(f_1,\mathbb{U}_1(\lambda,x))\\&=&
\mathbb{U}_1(\mathbb{U}_1(f_1,\lambda),x)\\&=&\mathbb{U}_1(a_1,x)
\end{eqnarray*}
for each $x\in [0,a_1]$.

Conversely, we prove the $(\alpha,\mathbb{U}_2)$-migrativity of $\mathbb{U}_1$ by the following four cases.

Case 1. If $x,y\in [0,a_1]$, then
$\mathbb{U}_1(\mathbb{U}_2(\alpha,x),y)=\mathbb{U}_1(\mathbb{U}_1(\lambda,x),y)=
\mathbb{U}_1(x,\mathbb{U}_1(\lambda,y))=\mathbb{U}_1(x,\mathbb{U}_2(\alpha,y))$.

Case 2. If $x,y\in [a_1,1]$, then
$\mathbb{U}_1(\mathbb{U}_2(\alpha,x),y)=\mathbb{U}_1(\mathbb{U}_1(\mu,x),y)=
\mathbb{U}_1(x,\mathbb{U}_1(\mu,y))=\mathbb{U}_1(x,\mathbb{U}_2(\alpha,y))$.

Case 3. If $x\in [0,a_1]$ and $y\in [a_1,1]$, then from Proposition \ref{pro-2.3} (iv) and $\mathbb{U}_1(a_1,x)=\mathbb{U}_1(\mu,x)$ for each $x\in [0,a_1]$, we have
\begin{eqnarray*}
\mathbb{U}_1(\mathbb{U}_2(\alpha,x),y)&=&\mathbb{U}_1(\mathbb{U}_1(\lambda,x),y)\\&=&
\mathbb{U}_1(x,\mathbb{U}_1(\lambda,y))\\&=&\mathbb{U}_1(x,\mathbb{U}_1(a_1,y))\\&=&\mathbb{U}_1(\mathbb{U}_1(x,a_1),y)
\\&=&\mathbb{U}_1(\mathbb{U}_1(x,\mu),y)\\&=&
\mathbb{U}_1(x,\mathbb{U}_1(\mu,y))\\&=&\mathbb{U}_1(x,\mathbb{U}_2(\alpha,y)).
\end{eqnarray*}

Case 4. If $x\in [a_1,1]$ and $y\in [0,a_1]$, then similar to Case 3, we obtain \eqref{eq-def}.
  \item
Assume that $\mathbb{U}_1$ is $(\alpha, \mathbb{U}_2)$-migrative. Then by Lemma \ref{lemma3.2} (ii) we get $\mathbb{U}_2(\alpha,a_1)\leq a_1$.
Since $\mathbb{U}_2(\alpha,x)\leq \mathbb{U}_2(\alpha,a_1)\leq a_1$ for any $x\leq a_1$, we have $\mathbb{U}_2(\alpha,x)=\mathbb{U}_1(x,\lambda)$ and $\lambda=\mathbb{U}_2(\alpha,e_1)\leq \mathbb{U}_2(\alpha,a_1)\leq a_1$. From Proposition \ref{pro-2.3} (i) we get $\mathbb{U}_2(\alpha,x)\geq \mathbb{U}_2(\alpha,a_1)=\mathbb{U}_1(\lambda,a_1)=\mathbb{U}_1(\lambda,f_1)=
\mathbb{U}_1(\mathbb{U}_2(\alpha,e_1),f_1)=\mathbb{U}_1(e_1,\mathbb{U}(\alpha,f_1))=a_1$ for any $x\geq a_1$. Consequently,
$$
\mathbb{U}_2(\alpha,x)=
\begin{cases}
\mathbb{U}_1(\lambda,x), & x\in [0,a_1],\\
\mathbb{U}_1(\mu,x), & x\in [a_1,1].\\
\end{cases}
$$

Next, from \eqref{eq-def} we have
\begin{eqnarray*}
\mathbb{U}_1(\lambda,x)&=&\mathbb{U}_1(\mathbb{U}_1(\lambda,e_1),x)\\&=&\mathbb{U}_1(\mathbb{U}_2(\alpha,e_1)
,x)\\&=&\mathbb{U}_1(e_1,\mathbb{U}_2(\alpha,x))\\&=&\mathbb{U}_1(e_1,\mathbb{U}_1(\mu,x))\\&=&
\mathbb{U}_1(\mathbb{U}_1(e_1,\mu),x)\\&=&\mathbb{U}_1(a_1,x)
\end{eqnarray*}
for each $x\in [a_1,1]$.

Conversely, we prove the $(\alpha,\mathbb{U}_2)$-migrativity of $\mathbb{U}_1$ by the following four cases.

Case 1. If $x,y\in [0,a_1]$, then
$\mathbb{U}_1(\mathbb{U}_2(\alpha,x),y)=\mathbb{U}_1(\mathbb{U}_1(\lambda,x),y)=
\mathbb{U}_1(x,\mathbb{U}_1(\lambda,y))=\mathbb{U}_1(x,\mathbb{U}_2(\alpha,y))$.

Case 2. If $x,y\in [a_1,1]$, then
$\mathbb{U}_1(\mathbb{U}_2(\alpha,x),y)=\mathbb{U}_1(\mathbb{U}_1(\mu,x),y)=
\mathbb{U}_1(x,\mathbb{U}_1(\mu,y))=\mathbb{U}_1(x,\mathbb{U}_2(\alpha,y))$.

Case 3. If $(x,y)\in [0,a_1]\times[a_1,1]$, then from Proposition \ref{pro-2.3} (i) and $\mathbb{U}_1(\lambda,x)=\mathbb{U}_1(a_1,x)$ for each $x\in [a_1,1]$,
\begin{eqnarray*}
\mathbb{U}_1(\mathbb{U}_2(\alpha,x),y)&=&\mathbb{U}_1(\mathbb{U}_1(\lambda,x),y)\\
&=&\mathbb{U}_1(x,\mathbb{U}_1(\lambda,y))\\&=&\mathbb{U}_1(x,\mathbb{U}_1(a_1,y))\\
&=&\mathbb{U}_1(\mathbb{U}_1(x,a_1),y)\\&=&\mathbb{U}_1(\mathbb{U}_1(x,\mu),y)\\&=&
\mathbb{U}_1(x,\mathbb{U}_1(\mu,y))\\&=&\mathbb{U}_1(x,\mathbb{U}_2(\alpha,y)).
\end{eqnarray*}

Case 4. If $(x,y)\in [a_1,1]\times [0,a_1]$, then the proof is in analogy to the Case 3.
\end{enumerate}
\end{proof}

Therefore, taking $\mathbb{U}_2(\alpha,e_1)$ as a classified point, we have the following theorem from Lemma \ref{lemma3.1}, Propositions \ref{pro-3.1} and \ref{pro-3.3}.
\begin{theorem}\label{thm3.1}
Let $\mathbb{U}_1$ and $\mathbb{U}_2$ be two 2-uninorms on the unit interval $[0,1]$ with 2-neutral elements $\{e_1,f_1\}_{a_1}$ and $\{e_2,f_2\}_{a_2}$, respectively.
\begin{enumerate}[{\rm (i)}]
  \item
If $\mathbb{U}_2(\alpha,e_1)\leq e_1$ and $\mathbb{U}_2(\alpha,f_1)\leq a_1$, then $\mathbb{U}_1$ is $(\alpha,\mathbb{U}_2)$-migrative iff $\mathbb{U}_1(\lambda,x)=\mathbb{U}_2(\alpha,x)$ for each $x\in [0,1]$.
  \item
If $\mathbb{U}_2(\alpha,e_1)\leq e_1$ and $\mathbb{U}_2(\alpha,f_1)\geq a_1$, then $\mathbb{U}_1$ is $(\alpha, \mathbb{U}_2)$-migrative iff
\begin{equation*}
\mathbb{U}_2(\alpha,x)=
\begin{cases}
\mathbb{U}_1(\lambda,x), & x\in [0,{a_1}],\\
\mathbb{U}_1(\mu,x), & x\in [{a_1},1]\\
\end{cases}
\end{equation*}
and
$\mathbb{U}_1(\lambda,x)=\mathbb{U}_1({a_1},x)$ for each $x\in [a_1,1]$.
  \item
If $\mathbb{U}_2(\alpha,e_1)\in [e_1,a_1]$, then $\mathbb{U}_1$ is $(\alpha, \mathbb{U}_2)$-migrative iff
\begin{equation*}
\mathbb{U}_2(\alpha,x)=
\begin{cases}
\mathbb{U}_1(\lambda,x), & x\in [0,{a_1}],\\
\mathbb{U}_1(\mu,x), & x\in [{a_1},1]\\
\end{cases}
\end{equation*}
and $\mathbb{U}_1(\mu,x)=\mathbb{U}_1({a_1},x)$ for each $x\in [0,a_1]$.
  \item
If $\mathbb{U}_2(\alpha,e_1)\geq a_1$, then $\mathbb{U}_1$ is $(\alpha,\mathbb{U}_2)$-migrative iff $\mathbb{U}_1(\mu,x)=\mathbb{U}_2(\alpha,x)$ for each $x\in [0,1]$.
\end{enumerate}
\end{theorem}

Taking $\mathbb{U}_2(\alpha,f_1)$ as a classified point, we have the following theorem from Lemma \ref{lemma3.1}, Proposition \ref{pro-3.1} and \ref{pro-3.3}.
\begin{theorem}\label{thm3.2}
Let $\mathbb{U}_1$ and $\mathbb{U}_2$ be two 2-uninorms on the unit interval $[0,1]$ with 2-neutral elements $\{e_1,f_1\}_{a_1}$ and $\{e_2,f_2\}_{a_2}$, respectively.
\begin{enumerate}[{\rm (i)}]
  \item
If $\mathbb{U}_2(\alpha,f_1)\leq a_1$, then $\mathbb{U}_1$ is $(\alpha,\mathbb{U}_2)$-migrative iff $\mathbb{U}_1(\lambda,x)=\mathbb{U}_2(\alpha,x)$ for any $x\in [0,1]$.
  \item
If $\mathbb{U}_2(\alpha,f_1)\in [a_1,f_1]$, then $\mathbb{U}_1$ is $(\alpha, \mathbb{U}_2)$-migrative iff
\begin{equation*}
\mathbb{U}_2(\alpha,x)=
\begin{cases}
\mathbb{U}_1(\lambda,x), & x\in [0,{a_1}],\\
\mathbb{U}_1(\mu,x), & x\in [{a_1},1]\\
\end{cases}
\end{equation*}
and
$\mathbb{U}_1(\lambda,x)=\mathbb{U}_1({a_1},x)$ for any $x\in [{a_1},1]$.
  \item
If $\mathbb{U}_2(\alpha,f_1)\geq f_1$ and $\mathbb{U}_2(\alpha,e_1)\leq a_1$, then $\mathbb{U}_1$ is $(\alpha, \mathbb{U}_2)$-migrative iff
\begin{equation*}
\mathbb{U}_2(\alpha,x)=
\begin{cases}
\mathbb{U}_1(\lambda,x), & x\in [0,{a_1}],\\
\mathbb{U}_1(\mu,x), & x\in [{a_1},1]\\
\end{cases}
\end{equation*}
and $\mathbb{U}_1(\mu,x)=\mathbb{U}_1({a_1},x)$ for any $x\in [0,{a_1}]$.
  \item
If $\mathbb{U}_2(\alpha,f_1)\geq f_1$ and $\mathbb{U}_2(\alpha,e_1)\geq a_1$, then $\mathbb{U}_1$ is $(\alpha,\mathbb{U}_2)$-migrative iff $\mathbb{U}_1(\mu,x)=\mathbb{U}_2(\alpha,x)$ for any $x\in [0,1]$.
\end{enumerate}
\end{theorem}

In particular, from Theorems \ref{thm3.1} and \ref{thm3.2}, we have the following remark.
\begin{remark}\label{remark3.1}
\begin{enumerate}[{\rm (i)}]
  \item
\emph{
Let $e_1=e_2=a_1=1$ in Theorem \ref{thm3.1}. Then both $\mathbb{U}_1$ and $\mathbb{U}_2$ are t-norms. Therefore, we get Theorem 3 of \cite{FR2011}.
  \item
Let $f_1=f_2=a_1=0$ in Theorem \ref{thm3.2}. Then both $\mathbb{U}_1$ and $\mathbb{U}_2$ are t-conorms. Thus $\mathbb{U}_1$ is $(\alpha, \mathbb{U}_2)$-migrative iff $\mathbb{U}_2(\alpha,x)=\mathbb{U}_1(\alpha,x)$ for each $x\in [0,1]$.
  \item
Let $e_2=a_1=1$ in Theorem \ref{thm3.2}. Then $\mathbb{U}_1$ is a uninorm and  $\mathbb{U}_2$ is a t-norms. Therefore, Theorem \ref{thm3.2} implies Theorem 3.1 of \cite{ZWY2021b}. In particular, if $\mathbb{U}_1$ is locally internal on the boundary (i.e., $\mathbb{U}_1(1,x)\in \{1,x\}$ for each $x\in [0,1]$ when $\mathbb{U}_1$ is a conjunctive uninorm \cite{MMRT2015}), then by $1>\alpha=\mathbb{U}_2(\alpha,1)=\mathbb{U}_1(\mathbb{U}_2(\alpha,1),e_1)=\mathbb{U}_1
(1,\mathbb{U}_2(\alpha,e_1))=\mathbb{U}_2(\alpha,e_1)$ for each $\alpha<1$ we obtain Proposition 5 of \cite{MMRT2015}.
  \item
Let $a_1=f_2=0$ in Theorem \ref{thm3.1}. Then $\mathbb{U}_1$ is a uninorm and  $\mathbb{U}_2$ is a t-conorm. Hence, we derive Theorem 3.3 of \cite{ZWY2021b}.
    \item
Let $a_1=a_2=1$ in Theorem \ref{thm3.2}. Then both $\mathbb{U}_1$ and $\mathbb{U}_2$ are uninorms and $\mathbb{U}_2(\alpha,1)\leq 1$. This deduces Lemma 1 in \cite{QR2015}.
  \item
Let $e_1=e_2=0$ and $f_1=f_2=1$ in Theorem \ref{thm3.1}. Then both $\mathbb{U}_1$ and $\mathbb{U}_2$ are nullnorms.
\begin{enumerate}[{\rm $\centerdot$}]
  \item
Let $a_1=a_2$. If $\alpha\leq a_1$, then $\mathbb{U}_2(\alpha,0)\leq a_1$, $\mathbb{U}_1(\mathbb{U}_2(\alpha,1),x)=\mathbb{U}_1(a_1,x)=a_1=
\mathbb{U}_1(\alpha,x)$ and $\mathbb{U}_1(\alpha,0)=\alpha$. If $\alpha>a_1$, then $\mathbb{U}_2(\alpha,0)=a_1$ and $\mathbb{U}_1(\alpha,1)=\alpha$. Consequently, we have Theorem 3 of \cite{ZSL2014}.
  \item
Let $a_1\neq a_2$. Assume that $a_1>a_2$ and $\mathbb{U}_1$ is $(\alpha,\mathbb{U}_2$)-migrative. Then from Theorem \ref{thm3.1} we can deduce that $\alpha\geq a_1$. Otherwise, $\alpha<a_1$. Then $\mathbb{U}_2(\alpha,0)<a_1$, which contradicts the fact $\mathbb{U}_2(\alpha,0)=\mathbb{U}_1(\mathbb{U}_2(\alpha,0),0)$ obtained from Theorem \ref{thm3.1} (iii). Assume that $a_1<a_2$ and $\mathbb{U}_1$ is $(\alpha,\mathbb{U}_2$)-migrative. Then from Theorem \ref{thm3.1} we can deduce that $\alpha\leq a_1$. Otherwise, $\alpha>a_1$. Then $\mathbb{U}_2(\alpha,0)=a_1$, contrary to the fact $\mathbb{U}_2(\alpha,1)=\mathbb{U}_1(\mathbb{U}_2(\alpha,0),1)$ obtained from Theorem \ref{thm3.1} (iv).
Thus, if both $\mathbb{U}_1$ and $\mathbb{U}_2$ are nullnorms with different absorbing elements, then we derive Theorem 7 in \cite{ZSL2014}.
\end{enumerate}
 \item
Let $e_1=0$, and $e_2=a_2=f_1=1$ in Theorem \ref{thm3.2}. Then $\mathbb{U}_1$ is a nullnorm and $\mathbb{U}_2$ is a t-norm. Since $\mathbb{U}_2(\alpha,1)\leq a_1$ yields $\alpha\leq a_1$ and $\mathbb{U}_2(\alpha,1)>a_1$ leads to $\alpha>a_1$, we obain the following three results:
\begin{enumerate}[{\rm $\centerdot$}]
  \item
If $\alpha\leq a_1$, then $\mathbb{U}_1$ is $(\alpha,\mathbb{U}_2$)-migrative iff $\mathbb{U}_2(\alpha,x)=\mathbb{U}_1(0,x)$ for any $x\in [0,1]$.
  \item
If $\alpha>a_1$, then since $\mathbb{U}_1(\lambda,y)=a_1=\mathbb{U}_1(a_1,y)$ for any $y\in [a_1,1]$, $\mathbb{U}_1$ is $(\alpha,\mathbb{U}_2$)-migrative iff
\begin{equation*}
\mathbb{U}_2(\alpha,x)
=
\begin{cases}
x,& x\in [0,a_1],\\
\mathbb{U}_1(\alpha,x), & x\in [a_1,1].\\
\end{cases}
\end{equation*}
  \item
In particular, if $\alpha=a_1$, then we have Theorem 3.5 in \cite{ZWY2021b}.
\end{enumerate}
  \item
If $e_1=f_2=0$ and $f_1=1$ in Theorem \ref{thm3.1}, then $\mathbb{U}_1$ is a nullnorm with absorbing element $a_1$ and $\mathbb{U}_2$ is a t-conorm. Since $\mathbb{U}_2(\alpha,0)< a_1$ leads to $\alpha< a_1$ and $\mathbb{U}_2(\alpha,0)\geq a_1$ yields $\alpha\geq a_1$, we have the following three results:
\begin{enumerate}[{\rm $\centerdot$}]
  \item
If $\alpha<a_1$, then from $\mathbb{U}_1(\mu,x)=a_1=\mathbb{U}_1(a_1,x)$ for any $x\in [0,a_1]$, $\mathbb{U}_1$ is $(\alpha,\mathbb{U}_2$)-migrative iff
\begin{equation*}
\mathbb{U}_2(\alpha,x)
=
\begin{cases}
\mathbb{U}_1(\alpha,x),& x\in [0,a_1],\\
x, & x\in [a_1,1].\\
\end{cases}
\end{equation*}
  \item
If $\alpha\geq a_1$, then $\mathbb{U}_1$ is $(\alpha,\mathbb{U}_2$)-migrative iff $\mathbb{U}_2(\alpha,x)=\mathbb{U}_1(1,x)$ for any $x\in [0,1]$.
\item
In particular, if $\alpha=a_1$, then we obtain Theorem 3.7 in \cite{ZWY2021b}.
\end{enumerate}
   \item
Let $f_1=f_2=1$ in Theorem \ref{thm3.2}. Then both $\mathbb{U}_1$ and $\mathbb{U}_2$ are uni-nullnorms with 2-neutral elements $\{e_1,1\}_{a_1}$ and $\{e_2,1\}_{a_2}$, respectively. Obviously, by Proposition \ref{pro-2.3} (i) we have $\mathbb{U}_1(\lambda,x)=\mathbb{U}_1(\lambda,a_1)=\mathbb{U}_1(\mu,a_1)=a_1=\mathbb{U}_1(a_1,x)$ for each $x\in [a_1,1]$ and $\mu>a_1$. Consequently, by Theorem \ref{thm3.2} (i) and (ii) we have all the results in  Section 4 of \cite{SQ2021}.
\item
Let $e_1=e_2=0$ in Theorem \ref{thm3.1}. Then both $\mathbb{U}_1$ and $\mathbb{U}_2$ are null-uninorms with 2-neutral elements $\{0,f_1\}_{a_1}$ and $\{0,f_2\}_{a_2}$, respectively.
Obviously, by Proposition \ref{pro-2.3} (iii) we have $\mathbb{U}_1(\mu,x)=\mathbb{U}_1(\mu,a_1)=\mathbb{U}_1(\lambda,a_1)
=a_1=\mathbb{U}_1(a_1,x)$ for any $x\in [0,a_1]$ and $\lambda<a_1$. Consequently, by Theorem \ref{thm3.1} (iii) and (iv)
we obtain all the results in Section 4 of \cite{SQ2021}.
\item
Let $a_1=a_2$ in Theorems \ref{thm3.1} and \ref{thm3.2}. Then $\mathbb{U}_1$ and $\mathbb{U}_2$ be two 2-uninorms on the unit interval $[0,1]$ with 2-neutral elements $\{e_1,f_1\}_{a_1}$ and $\{e_2,f_2\}_{a_1}$, respectively. If $\alpha\in [0,e_2]$, then $\lambda=\mathbb{U}_2(\alpha,e_1)\leq \mathbb{U}_2(e_2,e_1)=e_1\leq a_1$. If $\alpha\in [f_2,1]$, then $\mu=\mathbb{U}_2(\alpha,f_1)\geq \mathbb{U}_2(f_2,f_1)=f_1\geq a_1$. If $\alpha\in [e_2,f_2]$, then $\mathbb{U}_1(a_1,x)=\mathbb{U}_1(\lambda,x)$ for any $x\in [a_1,1]$ and $\mathbb{U}_1(a_1,x)=\mathbb{U}_1(\mu,x)$ for any $x\in [0,a_1]$. Therefore, we obtain Theorem 3.1 of \cite{LW2023}.
  \item
Let $e_1=e_2$ in Theorems \ref{thm3.1} and \ref{thm3.2}. Then $\mathbb{U}_1$ and $\mathbb{U}_2$ be two 2-uninorms on the unit interval $[0,1]$ with 2-neutral elements $\{e_1,f_1\}_{a_1}$ and $\{e_1,f_2\}_{a_2}$, respectively. Since $\lambda=\mathbb{U}_2(\alpha,e_1)=\alpha$ for any $\alpha\in [0,e_1]$, we have Theorems 3.2 and 3.3 in \cite{LW2023} in analogy to (xi).
  \item
Let $f_1=f_2$ in Theorems \ref{thm3.1} and \ref{thm3.2}. Then $\mathbb{U}_1$ and $\mathbb{U}_2$ be two 2-uninorms on the unit interval $[0,1]$ with 2-neutral elements $\{e_1,f_1\}_{a_1}$ and $\{e_2,f_1\}_{a_2}$, respectively. Since $\mu=\mathbb{U}_2(\alpha,f_1)=\alpha$ for any $\alpha\in [f_1,1]$, we get Theorems 3.4 and 3.5 of \cite{LW2023} in analogy to (xi).}
\end{enumerate}
\end{remark}

Further, we have the following corollary.
\begin{corollary}\label{coro3.1}
Let $\mathbb{U}$ be a 2-uninorm with 2-neutral elements $\{e_{\mathbb{U}},f_{\mathbb{U}}\}_{a_{\mathbb{U}}}$, $T$ be a t-norm, $S$ be a t-conorm, $U$ be a uninorm with a neutral element $e_{U}$, $N$ be a nullnorm with an absorbing element $a_N$, $F$ be a uni-nullnorm with 2-neutral elements $\{e_F,1\}_{a_F}$ and $G$ be a null-uninorm with 2-neutral elements $\{0,f_G\}_{a_G}$ on the unit interval $[0,1]$, respectively.
\begin{enumerate}[{\rm (i)}]
  \item
$T$ is $(\alpha,\mathbb{U})$-migrative iff $T(\mathbb{U}(\alpha,1),x)=\mathbb{U}(\alpha,x)$ for each $x\in [0,1]$..
  \item
\begin{enumerate}[{\rm $\centerdot$}]
  \item
If $T(\alpha,f_{\mathbb{U}})\leq a_{\mathbb{U}}$, then $\mathbb{U}$ is $(\alpha,T$)-migrative iff $\mathbb{U}(T(\alpha,e_{\mathbb{U}}),x)=T(\alpha,x)$ for each $x\in [0,1]$.
  \item
If $T(\alpha,f_{\mathbb{U}})\geq a_{\mathbb{U}}$, then $\mathbb{U}$ is $(\alpha,T$)-migrative iff
\begin{equation*}
T(\alpha,x)=
\begin{cases}
\mathbb{U}(T(\alpha,e_{\mathbb{U}}),x), & x\in [0,a_{\mathbb{U}}],\\
\mathbb{U}(T(\alpha,f_{\mathbb{U}}),x), & x\in [a_{\mathbb{U}},1],\\
\end{cases}
\end{equation*}
and
$\mathbb{U}(T(\alpha,e_{\mathbb{U}}),x)=\mathbb{U}(a_{\mathbb{U}},x)$ for any $x\in [a_{\mathbb{U}},1]$.
\end{enumerate}
  \item
$S$ is $(\alpha,\mathbb{U})$-migrative iff $S(\mathbb{U}(\alpha,0),x)=\mathbb{U}(\alpha,x)$ for each $x\in [0,1]$.
  \item
\begin{enumerate}[{\rm $\centerdot$}]
  \item
If $S(\alpha,e_{\mathbb{U}})\leq a_{\mathbb{U}}$, then $\mathbb{U}$ is $(\alpha,S$)-migrative iff \begin{equation*}
S(\alpha,x)=
\begin{cases}
\mathbb{U}(S(\alpha,e_{\mathbb{U}}),x), & x\in [0,a_{\mathbb{U}}],\\
\mathbb{U}(S(\alpha,f_{\mathbb{U}}),x), & x\in [a_{\mathbb{U}},1],\\
\end{cases}
\end{equation*}
and $\mathbb{U}(S(\alpha,f_{\mathbb{U}}),x)=\mathbb{U}(a_{\mathbb{U}},x)$ for any $x\in [0,a_{\mathbb{U}}]$.
  \item
If $S(\alpha,e_{\mathbb{U}})\geq a_{\mathbb{U}}$, then $\mathbb{U}$ is $(\alpha,S$)-migrative iff
$\mathbb{U}(S(\alpha,f_{\mathbb{U}}),x)=S(\alpha,x)$ for each $x\in [0,1]$.
\end{enumerate}
  \item
$U$ is $(\alpha,\mathbb{U})$-migrative iff $U(\mathbb{U}(\alpha,e_U),x)=\mathbb{U}(\alpha,x)$ for each $x\in [0,1]$.
  \item
\begin{enumerate}[{\rm $\centerdot$}]
  \item
If $\mathbb{U}(\alpha,0)\leq a_N$, then $N$ is $(\alpha,\mathbb{U}$)-migrative iff \begin{equation*}
\mathbb{U}(\alpha,x)=
\begin{cases}
N(\mathbb{U}(\alpha,0),x), & x\in [0,a_N],\\
N(\mathbb{U}(\alpha,1),x), & x\in [a_N,1].\\
\end{cases}
\end{equation*}
  \item
If $\mathbb{U}(\alpha,0)\geq a_N$, then $N$ is $(\alpha,\mathbb{U}$)-migrative iff
$N(\mathbb{U}(\alpha,1),x)=\mathbb{U}(\alpha,x)$ for each $x\in [0,1]$.
\end{enumerate}
  \item
\begin{enumerate}[{\rm $\centerdot$}]
  \item
If $\mathbb{U}(\alpha,1)\leq a_F$, then $F$ is $(\alpha,\mathbb{U}$)-migrative iff
$F(\mathbb{U}(\alpha,e_F),x)=\mathbb{U}(\alpha,x)$ for each $x\in [0,1]$.
  \item
If $\mathbb{U}(\alpha,1)\geq a_F$, then $F$ is $(\alpha,\mathbb{U}$)-migrative iff
\begin{equation*}
\mathbb{U}(\alpha,x)=
\begin{cases}
F(\mathbb{U}(\alpha,e_F),x), & x\in [0,a_F],\\
F(\mathbb{U}(\alpha,1),x), & x\in [a_F,1].\\
\end{cases}
\end{equation*}
\end{enumerate}
  \item
\begin{enumerate}[{\rm $\centerdot$}]
  \item
If $\mathbb{U}(\alpha,0)\leq a_G$, then $G$ is $(\alpha,\mathbb{U}$)-migrative iff
\begin{equation*}
\mathbb{U}(\alpha,x)=
\begin{cases}
G(\mathbb{U}(\alpha,0),x), & x\in [0,a_G],\\
G(\mathbb{U}(\alpha,f_G),x), & x\in [a_G,1].\\
\end{cases}
\end{equation*}
  \item
If $\mathbb{U}(\alpha,0)\geq a_G$, then $G$ is $(\alpha,\mathbb{U}$)-migrative iff
$G(\mathbb{U}(\alpha,f_G),x)=\mathbb{U}(\alpha,x)$ for each $x\in [0,1]$.
\end{enumerate}
\end{enumerate}
\end{corollary}
\begin{proof}
\begin{enumerate}[{\rm (i)}]
  \item
Taking $e_1=a_1=f_1=1$ in Theorem \ref{thm3.2}, then $\mathbb{U}_2(\alpha, 1)\leq 1$ for any $\alpha\leq 1$, and subsequently, from Theorem \ref{thm3.2} (i) we get that $\mathbb{U}_1$ is $(\alpha,\mathbb{U}_2)$-migrative iff $\mathbb{U}_1(\lambda,x)=\mathbb{U}_2(\alpha,x)$ for each $x\in [0,1]$. Finally, let $\mathbb{U}_1=T$ and $\mathbb{U}_2=\mathbb{U}$. Then we get (i).
  \item
Taking $e_2=a_2=f_2=1$ in Theorem \ref{thm3.1}, then $\mathbb{U}_2(\alpha,e_1)\leq e_1$ for any $\alpha\leq 1$, and subsequently, from Theorem \ref{thm3.1} (i) and (ii) we get that
\begin{enumerate}[{\rm $\centerdot$}]
  \item
If $\mathbb{U}_2(\alpha,f_1)\leq a_1$, then $\mathbb{U}_1$ is $(\alpha,\mathbb{U}_2$)-migrative iff $\mathbb{U}_1(\lambda,x)=\mathbb{U}_2(\alpha,x)$ for each $x\in [0,1]$.
  \item
If $\mathbb{U}_2(\alpha,f_1)\geq a_1$, then $\mathbb{U}_1$ is $(\alpha,\mathbb{U}_2$)-migrative iff
\begin{equation*}
\mathbb{U}_2(\alpha,x)=
\begin{cases}
\mathbb{U}_1(\lambda,x), & x\in [0,a_1],\\
\mathbb{U}_1(\mu,x), & x\in [a_1,1],\\
\end{cases}
\end{equation*}
and
$\mathbb{U}_1(\lambda,x)=\mathbb{U}_1({a_1},x)$ for any $x\in [{a_1},1]$.
\end{enumerate}
Finally, let $\mathbb{U}_1=\mathbb{U}$ and $\mathbb{U}_2=T$. Then (ii) holds.
  \item
Taking $e_1=a_1=f_1=0$ in Theorem \ref{thm3.1}, then $\mathbb{U}_2(\alpha,0)\geq 0$ for any $\alpha\geq 0$, and subsequently, from Theorem \ref{thm3.1} (iv) we get that $\mathbb{U}_1$ is $(\alpha,\mathbb{U}_2)$-migrative iff $\mathbb{U}_1(\mu,x)=\mathbb{U}_2(\alpha,x)$ for any $x\in [0,1]$. Finally, let $\mathbb{U}_1=S$ and $\mathbb{U}_2=\mathbb{U}$. Then we get (iii).
  \item
Taking $a_2=f_2=0$ in Theorem \ref{thm3.2}, then $\mathbb{U}_2(\alpha,f_1)\geq \mathbb{U}_2(0,f_1)=f_1$ for any $\alpha\geq 0$, and subsequently, from Theorem \ref{thm3.2} (iii) and (iv) we get that
\begin{enumerate}[{\rm $\centerdot$}]
  \item
If $\mathbb{U}_2(\alpha,e_1)\leq a_1$, then $\mathbb{U}_1$ is $(\alpha,\mathbb{U}_2$)-migrative iff
\begin{equation*}
\mathbb{U}_2(\alpha,x)=
\begin{cases}
\mathbb{U}_1(\lambda,x), & x\in [0,a_1],\\
\mathbb{U}_1(\mu,x), & x\in [a_1,1],\\
\end{cases}
\end{equation*}
and $\mathbb{U}_1(\mu,x)=\mathbb{U}_1({a_1},x)$ for any $x\in [0,{a_1}]$.
  \item
If $\mathbb{U}_2(\alpha,e_1)\geq a_1$, then $\mathbb{U}_1$ is $(\alpha,\mathbb{U}_2$)-migrative iff
$\mathbb{U}_1(\mu,x)=\mathbb{U}_2(\alpha,x)$ for each $x\in [0,1]$.
\end{enumerate}
Finally, let $\mathbb{U}_1=\mathbb{U}$ and $\mathbb{U}_2=S$. Then (iv) holds.
  \item
Taking $e_1=a_1=0$ in Theorem \ref{thm3.1}, then $\mathbb{U}_2(\alpha,0)\geq 0$ for any $\alpha\geq 0$, and subsequently, $\mathbb{U}_1$ is $(\alpha,\mathbb{U}_2)$-migrative iff $\mathbb{U}_1(\mu,x)=\mathbb{U}_2(\alpha,x)$ for any $x\in [0,1]$. Finally, let $\mathbb{U}_1=U$ and $\mathbb{U}_2=\mathbb{U}$. Then (v) holds.
  \item
Taking $e_1=0$ and $f_1=1$ in Theorem \ref{thm3.1}, then from Theorem \ref{thm3.1} (iii) and (iv) we get that
\begin{enumerate}[{\rm $\centerdot$}]
  \item
If $\mathbb{U}_2(\alpha,e_1)\leq a_1$, then $\mathbb{U}_1$ is $(\alpha,\mathbb{U}_2$)-migrative iff
\begin{equation*}
\mathbb{U}_2(\alpha,x)=
\begin{cases}
\mathbb{U}_1(\lambda,x), & x\in [0,a_1],\\
\mathbb{U}_1(\mu,x), & x\in [a_1,1].\\
\end{cases}
\end{equation*}
In this case, $\mathbb{U}_1(\mu,x)=\mathbb{U}_1({a_1},x)$ for each $x\in [0,a_1]$ is clearly established because of $\mathbb{U}_1(\mu,x)=\mathbb{U}_1(\mu,a_1)=\mathbb{U}_1(\lambda,a_1)=a_1=
\mathbb{U}_1({a_1},x)$.
  \item
If $\mathbb{U}_2(\alpha,e_1)\geq a_1$, then $\mathbb{U}_1$ is $(\alpha,\mathbb{U}_2$)-migrative iff
$\mathbb{U}_1(\mu,x)=\mathbb{U}_2(\alpha,x)$ for each $x\in [0,1]$.
\end{enumerate}
Finally, let $\mathbb{U}_1=N$ and $\mathbb{U}_2=\mathbb{U}$. Then (vi) holds.
   \item
Taking $f_1=1$ in Theorem \ref{thm3.2}, then from Theorem \ref{thm3.2} (i) and (ii) we get that
\begin{enumerate}[{\rm $\centerdot$}]
  \item
If $\mathbb{U}_2(\alpha,f_1)\leq a_1$, then $\mathbb{U}_1$ is $(\alpha,\mathbb{U}_2$)-migrative iff
$\mathbb{U}_1(\lambda,x)=\mathbb{U}_2(\alpha,x)$ for each $x\in [0,1]$.
  \item
If $\mathbb{U}_2(\alpha,f_1)\geq a_1$, then $\mathbb{U}_1$ is $(\alpha,\mathbb{U}_2$)-migrative iff
\begin{equation*}
\mathbb{U}_2(\alpha,x)=
\begin{cases}
\mathbb{U}_1(\lambda,x), & x\in [0,a_1],\\
\mathbb{U}_1(\mu,x), & x\in [a_1,1].\\
\end{cases}
\end{equation*}
In this case, $\mathbb{U}_1(\lambda,x)=\mathbb{U}_1({a_1},x)$ for any $x\in [{a_1},1]$ is clearly established because of $\mathbb{U}_1(\lambda,x)=\mathbb{U}_1(\lambda,a_1)=\mathbb{U}_1(\mu,a_1)=a_1=\mathbb{U}_1({a_1},x)$.
\end{enumerate}
Finally, let $\mathbb{U}_1=F$ and $\mathbb{U}_2=\mathbb{U}$. Then (vii) holds.
   \item
Taking $e_1=0$ in Theorem \ref{thm3.1}, then from Theorem \ref{thm3.1} (iii) and (iv) we get that
\begin{enumerate}[{\rm $\centerdot$}]
  \item
If $\mathbb{U}_2(\alpha,e_1)\leq a_1$, then $\mathbb{U}_1$ is $(\alpha,\mathbb{U}_2$)-migrative iff
\begin{equation*}
\mathbb{U}_2(\alpha,x)=
\begin{cases}
\mathbb{U}_1(\lambda,x), & x\in [0,a_1],\\
\mathbb{U}_1(\mu,x), & x\in [a_1,1].\\
\end{cases}
\end{equation*}
In this case, $\mathbb{U}_1(\mu,x)=\mathbb{U}_1({a_1},x)$ for each $x\in [0,a_1]$ is clearly established because of $\mathbb{U}_1(\mu,x)=\mathbb{U}_1(\mu,a_1)=\mathbb{U}_1(\lambda,a_1)=a_1=\mathbb{U}_1({a_1},x)$.
  \item
If $\mathbb{U}_2(\alpha,e_1)\geq a_1$, then $\mathbb{U}_1$ is $(\alpha,\mathbb{U}_2$)-migrative iff
$\mathbb{U}_1(\mu,x)=\mathbb{U}_2(\alpha,x)$ for each $x\in [0,1]$.
\end{enumerate}
Finally, let $\mathbb{U}_1=G$ and $\mathbb{U}_2=\mathbb{U}$. Then (viii) holds.
\end{enumerate}
\end{proof}

As a conclusion of this section, we give a simple example to help us understand the above characterizations.
\begin{example}
\emph{
Define $\mathbb{U}_1$ and $\mathbb{U}_2$ as follows:
\begin{equation*}
\mathbb{U}_1(x,y)=
\begin{cases}
\min(x,y),&  (x,y)\in [0,e_1]^2\cup [a_1,f_1]^2\\
a_1,& (x,y)\in [0,a_1]\times [a_1,f_1]\cup [a_1,f_1]\times [0,a_1],\\
\max(x,y), & \rm{otherwise}\\
\end{cases}
\end{equation*}
and
\begin{equation*}
\mathbb{U}_2(x,y)=
\begin{cases}
\min(x,y),&  (x,y)\in [0,e_2]^2\cup [a_2,f_2]^2\\
a_2,& (x,y)\in  [0,a_2]\times [a_2,f_2]\cup [a_2,f_2]\times [0,a_2],\\
\max(x,y), & \rm{otherwise}.\\
\end{cases}
\end{equation*}
One can check that $\mathbb{U}_1$ and $\mathbb{U}_2$ are 2-uninorms with 2-neutral elements $\{e_1,f_1\}_{a_1}$ and $\{e_2,f_2\}_{a_2}$, respectively. Let $e_1=0.2$, $a_1=0.6$, $f_1=0.8$, $e_2=0.3$, $a_2=0.5$ and $f_2=0.7$.}
\begin{enumerate}[{\rm $\diamond$}]
  \item
\emph{If $\alpha\in [0,0.3]$, then $\mathbb{U}_2(\alpha,e_1)\leq e_1$, $\mu=\mathbb{U}_2(\alpha,f_1)=0.8= f_1$ and $\mathbb{U}_2(\alpha,0.6)=0.5\neq 0.6=\mathbb{U}_1(\mu,0.6)$, showing that $\mathbb{U}_1$ is not $(\alpha,\mathbb{U}_2$)-migrative by Proposition \ref{pro-3.3}.
  \item
If $\alpha\in (0.3,0.5]$, then $\mathbb{U}_2(\alpha,e_1)=\alpha\in [e_1,a_1]$, $\mu=\mathbb{U}_2(\alpha,f_1)=0.8$ and $\mathbb{U}_2(\alpha,0.6)=0.5\neq 0.6=\mathbb{U}_1(\mu,0.6)$, showing that $\mathbb{U}_1$ is not $(\alpha,\mathbb{U}_2$)-migrative by Proposition \ref{pro-3.3}.
  \item
If $\alpha\in (0.5,0.7)$, then $\mathbb{U}_2(\alpha,e_1)=0.5$, $\mu=\mathbb{U}_2(\alpha,f_1)=0.8=f_1$ and $\mathbb{U}_2(\alpha,0.7)=\alpha\neq 0.7=\mathbb{U}_1(\mu,0.7)$, showing that $\mathbb{U}_1$ is not $(\alpha,\mathbb{U}_2$)-migrative by Proposition \ref{pro-3.3}.
  \item
If $\alpha=0.7$, then $\mathbb{U}_2(\alpha,e_1)=0.5$, $\mu=\mathbb{U}_2(\alpha,f_1)=0.8=f_1$, $\mathbb{U}_2(\alpha,x)=0.5=\mathbb{U}_1(\lambda,x)$ for each $x\in [0,0.5]$, $\mathbb{U}_2(\alpha,x)=x=\mathbb{U}_1(\lambda,x)$ for each $x\in [0.5,0.6]$ and $\mathbb{U}_2(\alpha,x)=x=\mathbb{U}_1(\mu,x)$ for each $x\in [0.6,1]$, showing that $\mathbb{U}_1$ is $(\alpha,\mathbb{U}_2)$-migrative by Proposition \ref{pro-3.3}.
  \item
If $\alpha\in (0.7,0.8]$, then $\mathbb{U}_2(\alpha,e_1)=\alpha\geq a_1$, $\mu=\mathbb{U}_2(\alpha,f_1)=0.8$ and $\mathbb{U}_2(\alpha,0.7)=\alpha\neq 0.7=\mathbb{U}_1(\mu,0.7)$, showing that $\mathbb{U}_1$ is not $(\alpha,\mathbb{U}_2$)-migrative by Proposition \ref{pro-3.1}.
  \item
If $\alpha\in (0.8,1]$, then $\mathbb{U}_2(\alpha,e_1)=\alpha\geq a_1$, $\mu=\mathbb{U}_2(\alpha,f_1)=\alpha$ and $\mathbb{U}_2(\alpha,x)=\alpha\vee x=\mathbb{U}_1(\mu,x)$ for each $x\in [0,1]$, showing that $\mathbb{U}_1$ is $(\alpha,\mathbb{U}_2$)-migrative by Proposition \ref{pro-3.1}.}
\end{enumerate}
\end{example}

\section{Conclusions}\label{sec6}
This article provided a complete characterization for the $\alpha$-migrativity of 2-uninorms with $\alpha\in [0,1]$. Comparing with the results of Liang and Wang \cite{LW2023}, we supplemented the migrativity of 2-uninorms $\mathbb{U}_1$ and $\mathbb{U}_2$ when $e_1\neq e_2$ and $a_1\neq a_2$ and $f_1\neq f_2$ where $\{e_1,f_1\}_{a_1}$ and $\{e_2,f_2\}_{a_2}$ are 2-neutral elements of $\mathbb{U}_1$ and $\mathbb{U}_2$, respectively. It is worth pointing out that our results are suitable for all 2-uninorms, unlike the results of Huang and Qin \cite{HQ2021} and Wang et al. \cite{WZZL2019} which are for some special 2-uninorms. Furthermore, Corollary \ref{coro3.1} presented some characterizations of the migrativity between 2-uninorm and other subclasses of 2-uninorms, which together Remark \ref{remark3.1} shows us that our results unify and generalize some current existing characterizations for the $\alpha$-migrativity of t-norms, t-conorms, uninorms, nullnorms, uni-nullnorms and null-uninorms, respectively.

\section*{Declarations}
\section*{Conflict of interest}
No potential conflict of interest was declared by the authors. All data generated or analysed during this study are included in this
published article (and its supplementary information files).

\end{document}